\documentstyle[11pt,amssymb,graphicx]{article}
\begin{document}  
\date{} 
\renewcommand{\thefootnote}{}

\def\thebibliography#1{\begin{center}{\bf\normalsize References}\end{center}
 \list{[{\bf \arabic{enumi}}]}{\settowidth\labelwidth{[#1]}
 \leftmargin\labelwidth
 \advance\leftmargin\labelsep
 \usecounter{enumi}}
 \def\newblock{\hskip .11em plus .33em minus .07em}
 \sloppy\clubpenalty4000\widowpenalty4000
 \sfcode`\.=1000\relax}

\def\arf{{\mathrm Arf}}
\def\rarf{\overline{{\mathrm Arf}}}
\def\lk{{\mathrm lk}}

\title{\vspace*{-1.5cm}{\large CLASSIFICATION OF LINKS UP TO SELF $\#$-MOVE}}

\author{{\normalsize TETSUO SHIBUYA}\\
{\small Department of Mathematics, Osaka Institute of Technology}\\[-1mm]
{\small Omiya 5-16-1, Asahi, Osaka 535-8585, Japan}\\[-1mm]
{\small e-mail: shibuya@ge.oit.ac.jp}\\[3mm]
{\normalsize AKIRA YASUHARA }\\
{\small Department of Mathematics, Tokyo Gakugei University}\\[-1mm]
{\small Nukuikita 4-1-1, Koganei, Tokyo 184-8501, Japan}\\
{\small {\em Current address}: }\\[-1mm]
{\small Department of Mathematics, The George Washington University}\\[-1mm]
{\small Washington, DC 20052, USA}\\[-1mm]
{\small e-mail: yasuhara@u-gakugei.ac.jp}\\
}

\maketitle


\baselineskip=12pt
\vspace*{-3em}
\noindent
\begin{center}
{\small{\em Dedicated to Professor Shin'ich Suzuki 
for his 60th birthday}}
\end{center}

{\small 
\begin{quote}
\begin{center}A{\sc bstract}\end{center}
A pass-move and a $\#$-move are local moves on oriented links defined 
by L.H. Kauffman and H. Murakami respectively. 
Two links are self pass-equivalent (resp. self $\#$-equivalent) if 
one can be deformed into the other by pass-moves (resp. $\#$-moves), 
where non of them can occur between distinct components of the link.
These relations are equivalence relations on ordered oriented links and 
stronger than link-homotopy defined by J. Milnor. 
We give two complete classifications of links with 
arbitrarily many components
up to self pass-equivalence and up to self $\#$-equivalence respectively. 
So our classifications give subdivisions of link-homotopy classes.
\end{quote}}

\baselineskip=16pt

\begin{center}{\bf{1. Introduction}}\end{center}

We shall work in piecewise linear category. 
All links will be assumed 
to be ordered and oriented. 

A {\em pass-move} \cite{Kau} (resp. {\em $\#$-move} \cite{Mur}) 
is a local move on oriented links 
as illustrated in Figure 1.1(a) (resp. 1.1(b)). 
If the four strands in Figure 1.1(a) (resp. 1.1(b)) belong to 
the same component of a link, we call it a 
{\em self pass-move} (resp. {\em self $\#$-move}) 
\cite{C-F}, \cite{Shi2}, \cite{Shi3}. 
We note that the first author called pass-move and 
$\#$-move $\#$(II)-move and $\#$(I)-move respectively 
in his prior papers \cite{Shi2}, \cite{Shi3}, \cite{Shi4}, etc. 
Two links are {\em self pass-equivalent} 
(resp. {\em self $\#$-equivalent}) 
if one can be deformed into the other by 
a finite sequence of self pass-moves (resp. self $\#$-moves). 
Two links are {\em link-homotopic} if one can be deformed 
into the other by finite sequence of {\em self crossing changes} \cite{Mil}. 
Since both self pass-move and self $\#$-move are realized by 
self crossing changes, self pass-equivalence and self $\#$-equivalence 
are stronger than link-homotopy. Link-homotopy classification is 
already done by N. Habegger and X.S. Lin \cite{H-L}. 
In this paper we give two complete classifications of links 
with arbitrarily many components 
up to self pass-equivalence and up to self $\#$-equivalence respectively. 
So our classifications give subdivisions of link-homotopy classes.

\medskip
\begin{center} 
\includegraphics[trim=0mm 0mm 0mm 0mm, width=.9\linewidth]
{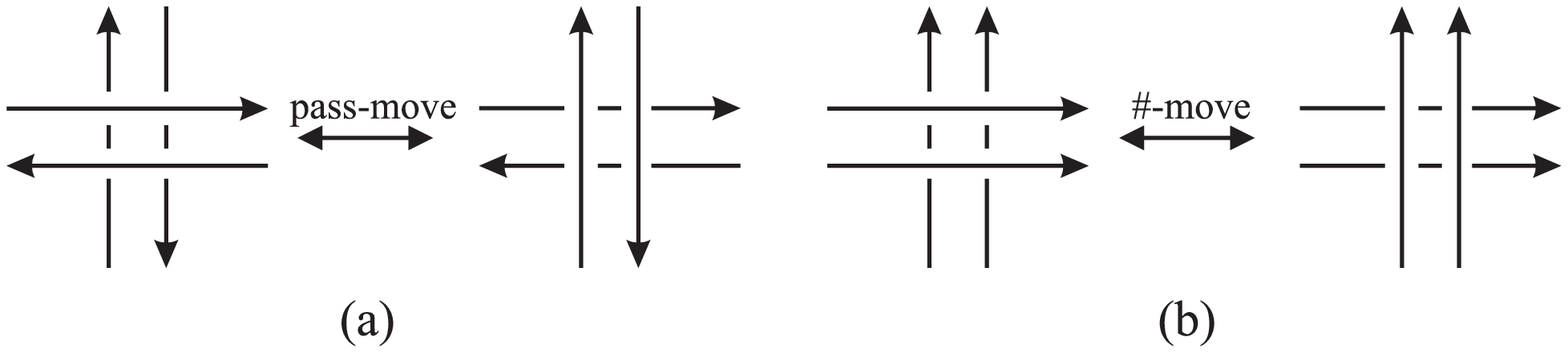}

Figure 1.1
\end{center}

An $n$-component link $l=k_1\cup\cdots\cup k_n$ is {\em proper} if 
the linking number $\lk(l-k_i,k_i)$ is even for 
any $i(=1,...,n)$. We define that a knot is a proper link. 
For a proper link $l=k_1\cup\cdots\cup k_n$, 
we call $\arf(l)-\sum_{i=1}^n\arf(k_i)$ (mod 2) 
the {\em reduced Arf invariant} \cite{Shi2} and denote it by 
$\rarf(l)$, where $\arf$ is the {\em Arf invariant} \cite{Rob}.

\medskip
{\bf Theorem 1.1.} {\em Let $l=k_1\cup\cdots\cup k_n$ and 
$l'=k'_1\cup\cdots\cup k'_n$ be $n$-component links. Then the 
following {\rm (i)} and {\rm (ii)} hold.
\begin{enumerate}
\item[{\rm (i)}] $l$ and $l'$ are self pass-equivalent if and only if 
they are link-homotopic and 
$\arf(k_{i_1}\cup\cdots\cup k_{i_p})=\arf(k'_{i_1}\cup\cdots\cup k'_{i_p})$ 
for any proper links $k_{i_1}\cup\cdots\cup k_{i_p}\subseteq l$ and 
$k'_{i_1}\cup\cdots\cup k'_{i_p}\subseteq l'$.
\item[{\rm (ii)}] $l$ and $l'$ are self $\#$-equivalent if and only if 
they are link-homotopic and 
$\rarf(k_{i_1}\cup\cdots\cup k_{i_p})=\rarf(k'_{i_1}\cup\cdots\cup k'_{i_p})$ 
for any proper links $k_{i_1}\cup\cdots\cup k_{i_p}\subseteq l$ and 
$k'_{i_1}\cup\cdots\cup k'_{i_p}\subseteq l'$.
\end{enumerate}
}

\medskip
For two-component links, both self pass-equivalence classification 
and self $\#$-equivalence classification are done by the first author 
\cite{Shi3}. His proof can be applied to only two-component links. 
So we need different approach to proving Theorem 1.1.

A link $l=k_1\cup\cdots\cup k_n$ is {\em ${\Bbb Z}_2$-algebraically split} 
if $\lk(k_i,k_j)$ is even for any $i,j\ (1\leq i< j\leq n)$.
We note that if $l=k_1\cup\cdots\cup k_n$ is 
${\Bbb Z}_2$-algebraically split link, then $l$ and 
$k_i\cup k_j$ $(1\leq i<j\leq n)$ are proper.

\medskip
{\bf Theorem 1.2.} {\em Let $l=k_1\cup\cdots\cup k_n$ and 
$l'=k'_1\cup\cdots\cup k'_n$ be $n$-component ${\Bbb Z}_2$-algebraically 
split links. If $l$ and $l'$ are link-homotopic, then 
\[\rarf(l)+\sum_{1\leq i<j\leq n}\rarf(k_i\cup k_j)
\equiv
\rarf(l')+
\sum_{1\leq i<j\leq n}\rarf(k'_i\cup k'_j)
\ (\mbox{mod $2$}).\]}

\bigskip
\begin{center}{\bf{2. Preliminaries}}\end{center}

In this section, we collect several results 
in order to prove Theorems 1.1 and 1.2.

Let $l=k_1\cup\cdots\cup k_n$ and 
$l'=k'_1\cup\cdots\cup k'_n$ be $n$-component links. 
Suppose that there is a disjoint union ${\cal A}=
A_1\cup \cdots\cup A_n$ of $n$ annuli in $S^3\times[0,1]$
with $(\partial (S^3\times[0,1]), \partial A_i)=
(S^3\times\{0\},k_i)\cup(-S^3\times\{1\},-k'_i)$ $(i=1,...,n)$
such that 
\begin{enumerate}
\item[(i)] ${\cal A}$ is locally flat except for finite points 
$p_1,...,p_m$ in the interior of $\cal A$, and 
\item[(ii)] for each $p_j$ ($j=1,2,...,m$), there is a small neighborhood 
$N(p_j)$ of $p_j$ in $S^3\times[0,1]$ such that 
$(\partial N(p_j), \partial(N(p_j)\cap {\cal A}))$ is a link 
as illustrated in Figure 2.1, 
\end{enumerate}
where $-X$ denotes $X$ with the opposite orientation.
Then $\cal A$ is called 
a {\em pass-annuli} between $l$ and $l'$.

\medskip
\begin{center} 
\includegraphics[trim=0mm 0mm 0mm 0mm, width=.3\linewidth]
{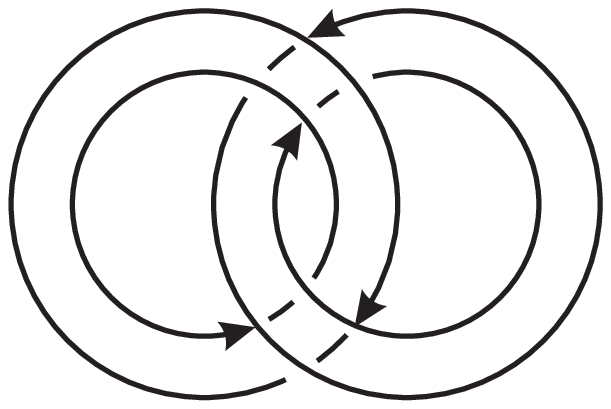}

Figure 2.1
\end{center}

The following is proved by the first author in \cite{Shi2}.

\medskip
{\bf Lemma 2.1.} {\em Two links $l$ and $l'$ are self pass-equivalent 
if and only if there is a pass-annuli between them. $\Box$}

\medskip
It is known that a pass-move is realized by a finite sequence 
of $\#$-moves \cite{M-N}. Thus we have the following.

\medskip
{\bf Lemma 2.2.} {\em If two links $l$ and $l'$ are self pass-equivalent, 
then they are self $\#$-equivalent. $\Box$}

\medskip
A {\em $\Gamma$-move} \cite{Kau} is a local move on oriented links 
as illustrated in Figure 2.2. 

\medskip
\begin{center} 
\includegraphics[trim=0mm 0mm 0mm 0mm, width=.75\linewidth]
{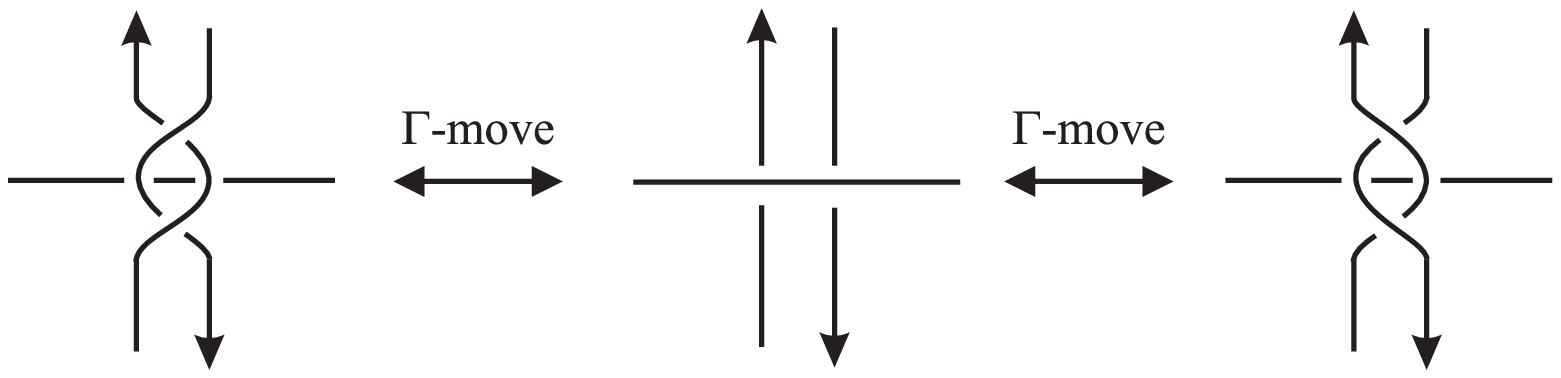}

Figure 2.2
\end{center}

The following is known \cite{Kau}.

\medskip
{\bf Lemma 2.3.} {\em A $\Gamma$-move is realized by a 
single pass-move. $\Box$}

\medskip
Let $l=k_1\cup\cdots\cup k_n$ and $l'=k'_1\cup\cdots\cup k'_n$ be 
$n$-component links such that there is a 3-ball $B^3$ in $S^3$ with 
$B^3\cap(l\cup l')=l$. Let $b_1,...,b_n$ be mutually disjoint disks in
$S^3$ such that $b_i\cap l=\partial b_i\cap k_i$ and 
$b_i\cap l'=\partial b_i\cap k'_i$ are arcs for each $i$.
Then the link $l\cup l'\cup(\bigcup_{i=1}^n \partial b_i) 
-(\bigcup {\mathrm int}(b_i\cap(l\cup l')))$ is called 
a {\em band sum} (or a {\em product fusion} \cite{Shi1}) 
of $l$ and $l'$ and denoted by 
$(k_1\#_{b_1}k'_1)\cup\cdots\cup(k_n\#_{b_n}k'_n)$. 
Note that a band sum of $l$ and $l'$ is ${\Bbb Z}_2$-algebraically 
split if $\lk(k_i,k_j)\equiv\lk(k'_i,k'_j)$ (mod 2) $(1\leq i<j\leq n)$. 

The following is proved by the first author in \cite{Shi1}. 

\medskip
{\bf Lemma 2.4.} {\em Two links $l$ and $l'$ are link-homotopic 
if and only if there is a band sum of $l$ and $-\overline{l'}$ 
that is link-homotopic to a trivial link, where 
$(S^3,-\overline{l'})\cong(-S^3,-l')$.  $\Box$}

\medskip
By the definition of the Arf invariant 
via 4-dimensional topology \cite{Rob}, we have the following.

\medskip
{\bf Lemma 2.5.} {\em Let $l$ and $l'$ be proper links and 
$L$ a band sum of $l$ and $-\overline{l'}$. 
Then $L$ is proper and $\arf(L)\equiv\arf(l)+\arf(l')$ 
$($mod $2$$)$. $\Box$ }

\medskip
The following lemma forms an interesting contrast to the lemma above. 

\medskip
{\bf Lemma 2.6.} {\em Let $l=k_1\cup k_2$ and $l'=k'_1\cup k'_2$ 
be $2$-component links with $\lk(k_1,k_2)$ and $\lk(k'_1,k'_2)$ odd. 
Let $L=(k_1\#_{b_1}(-\overline{k'_1}))\cup
(k_2\#_{b_2}(-\overline{k'_2}))$ be a band sum and $L'$ a band sum 
obtained from $L$ by adding a single full-twist to $b_2;$ see Figure $2.3$. 
Then $L$ and $L'$ are proper and link-homotopic, and 
$\arf(L)\neq\arf(L')$. }

\medskip
\begin{center} 
\includegraphics[trim=0mm 0mm 0mm 0mm, width=.5\linewidth]
{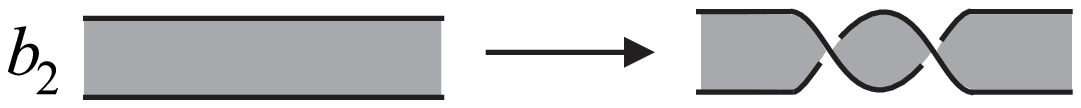}

Figure 2.3
\end{center}

\medskip
{\bf Proof.} 
Clearly $L$ and $L'$ are proper and link-homotopic. 
So we shall show $\arf(L)\neq\arf(L')$.

Let $a_i$ be the $i$th coefficient of the Conway polynomial. 
Then we have 
\[a_3(L)-a_3(L')=a_2((k_1\#_{b_1}(-\overline{k'_1}))
\cup k_2\cup(-\overline{k'_2})).\]
It is known that the third coefficient of the Conway polynomial 
of two-component proper link is mod 2 congruent to the sum of 
the Arf invariants of the link and the components 
\cite{Murasugi}. 
This and Lemma 2.5 imply 
$\arf(L)-\arf(L')\equiv a_3(L)-a_3(L')$ (mod 2). 
By \cite{Hos}, 
\[\begin{array}{l}
a_2((k_1\#_{b_1}(-\overline{k'_1}))
\cup k_2\cup(-\overline{k'_2}))\\
\hspace*{1cm}=
\lk(k_1\#_{b_1}(-\overline{k'_1}),k_2)
\lk(k_2,-\overline{k'_2})+
\lk(k_2,-\overline{k'_2})
\lk(-\overline{k'_2},k_1\#_{b_1}(-\overline{k'_1}))\\
\hspace*{7cm}+
\lk(-\overline{k'_2},k_1\#_{b_1}(-\overline{k'_1}))
\lk(k_1\#_{b_1}(-\overline{k'_1}),k_2).
\end{array}\]
Thus we have $\arf(L)-\arf(L')\equiv 1$ (mod 2). $\Box$

\medskip
A $\Delta$-move \cite{M-N} is a local move on links as 
illustrated in Figure 2.4. 
If at least two of the three strands in Figure 2.4 belong to the 
same component of a link, we call it a {\em quasi self $\Delta$-move} 
\cite{N-S2}. 
Two links are {\em quasi self $\Delta$-equivalent} if one can be 
deformed into the other by a finite sequence of quasi self $\Delta$-moves. 

\medskip
\begin{center} 
\includegraphics[trim=0mm 0mm 0mm 0mm, width=.45\linewidth]
{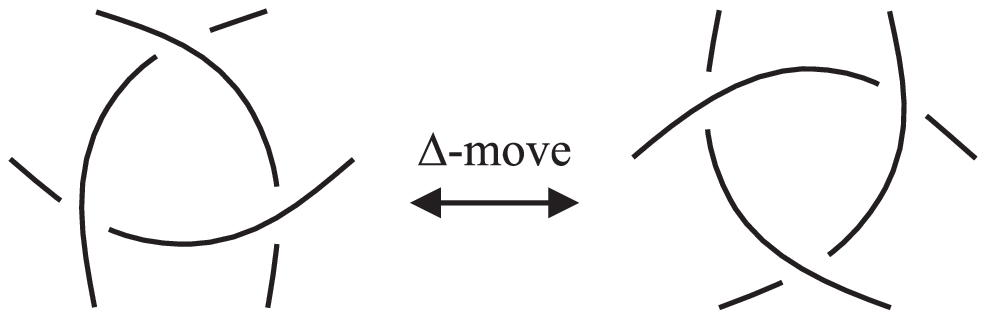}

Figure 2.4
\end{center}

The following is proved by Y. Nakanishi and the first author in \cite{N-S2}.

\medskip
{\bf Lemma 2.7.} {\em Two links are link-homotopic if and only if 
they are quasi self $\Delta$-equivalent. $\Box$}

\bigskip
\begin{center}{\bf{3. Proofs of Theorems 1.1 and 1.2}}\end{center}

{\bf Proof of Theorem 1.2.} 
Since $l$ is link-homotopic to $l'$, by Lemma 2.7, 
$l$ is quasi self $\Delta$-equivalent to $l'$.  
It is sufficient to consider the case that 
$l'$ is obtained from $l$ by a single 
quasi self $\Delta$-move. 

Suppose that the three strands of the $\Delta$-move 
that is applied to the deformation from $l$ into $l'$ 
belong to one component of $l$. 
Without loss of generality we may assume that the component is $k_1$. 
Note that $k_i$ and $k'_i$ are ambient isotopic for any $i(\neq1)$, 
and that $k_i\cup k_j$ and $k'_i\cup k'_j$ are ambient isotopic 
for any $i<j$ $(i\neq1)$. 
Since a $\Delta$-move changes the value of Arf invariant \cite{M-N}, 
we have $\arf(l)\neq\arf(l')$, $\arf(k_1)\neq\arf(k'_1)$ 
and $\arf(k_1\cup k_j)\neq\arf(k'_1\cup k'_j)$. 
Thus we have $\rarf(l)=\rarf(l')$ and 
$\rarf(k_1\cup k_j)=\rarf(k'_1\cup k'_j)$. 
So we have the conclusion.

We consider the other case, i.e., the three strands of the $\Delta$-move 
belong to exactly two components of $l$. 
Without loss of generality we 
may assume that the two components are $k_1$ and $k_2$. 
Note that $k_i$ and $k'_i$ are ambient isotopic for any $i$, 
and that $k_i\cup k_j$ and $k'_i\cup k'_j$ are ambient isotopic 
for any $i<j$ $((i,j)\neq(1,2))$. 
Since $\arf(l)\neq\arf(l')$ and $\arf(k_1\cup k_2)\neq\arf(k'_1\cup k'_2)$, 
$\arf(l)+\arf(k_1\cup k_2)\equiv\arf(l')+\arf(k'_1\cup k'_2)$ (mod 2). 
This completes the proof. $\Box$

\medskip
{\bf Lemma 3.1.} {\em Let $l=k_1\cup\cdots\cup k_n$ and 
$l'=k'_1\cup\cdots\cup k'_n$ be $n$-component 
${\Bbb Z}_2$-algebraically split links. If 
$l$ and $l'$ are link-homotopic, $\arf(k_i)=\arf(k'_i)$ 
$(i=1,...,n)$ and $\arf(k_i\cup k_j)=\arf(k'_i\cup k'_j)$ 
$(1\leq i<j\leq n)$, then $l$ and $l'$ are self pass-equivalent. }

\medskip
{\bf Proof.} 
Since $l$ is link-homotopic to $l'$, 
by Lemma 2.7, $l$ is quasi self $\Delta$-equivalent to $l'$. 
Let $u$ be the minimum number of quasi self $\Delta$-moves 
which are needed to deform $l$ into $l'$. 
By Theorem 1.2, $\arf(l)=\arf(l')$. 
Since a $\Delta$-move changes the value of the Arf invariant, 
$u$ is even. 
It is sufficient to consider the case $u=2$. 
Therefore there is a union ${\cal A}=A_1\cup\cdots\cup A_n$ 
of {\em level-preserving} $n$ annuli in $S^3\times[0,1]$ 
with $(\partial (S^3\times[0,1]), \partial A_i)=
(S^3\times\{0\},k_i)\cup(-S^3\times\{1\},-k'_i)$ $(i=1,...,n)$
such that 
\begin{enumerate}
\item[(i)] ${\cal A}$ is locally flat except for exactly two points 
$p_1,p_2$ in the interior of $\cal A$, and 
\item[(ii)] for each $p_t$ ($t=1,2$) there is a small neighborhood 
$N(p_t)$ of $p_t$ in $S^3\times[0,1]$ such that 
$(\partial N(p_t), \partial(N(p_t)\cap {\cal A}))$ is 
the Borromean ring $R_t$, at least two components of which 
belong to some $A_i$. 
\end{enumerate} 
A singular points $p_t$ is called {\em type $(i)$} 
if the three components of $R_t$ belong to $A_i$ and 
{\em type $(i,j)$} $(i<j)$ if one or two componets are in $A_i$ 
and the others in $A_j$. 
For each $i$ (resp. $i,j$), let $u_i$ (resp. $u_{i,j}$)
be the number of the singular points 
of type $(i)$ (resp. type $(i,j)$). 
We note that a number of $\Delta$-moves which are needed to 
deform $k_i$ into $k'_i$ (resp. $k_i\cup k_j$ into $k'_i\cup k'_j$) 
is equal to $u_i$ (resp. $u_{i,j}+u_i+u_j$). 
By the hypothesis of this lemma, we have $u_i$ and 
$u_{i,j}+u_i+u_j$ are even. 
Hence $u_i$ and $u_{i,j}$ are even. 
This implies that both $p_1$ and $p_2$ are the same type. 

Suppose that $p_1$ and $p_2$ are type $(i,j)$. 
Without loss of generality we may assume that $(i,j)=(1,2)$ and two 
components of the Borromean ring $R_1$ belong to $A_2$. 
Let $\alpha$ be an arc in the interior of $A_1$ that connects 
two singular points $p_1$ and $p_2$ of type $(1,2)$, and 
let $(S^3,L)=(\partial N(\alpha),
\partial(N(\alpha)\cap(A_1\cup A_2)))$. 
Then $L$ is a $5$-component link as illustrated in either 
Figure 3.1(a) or (b). 
In the case that $L$ is as Figure 3.1(a), we can deform 
$L$ into a trivial link by applying $\Gamma$-moves to the sublink 
$L\cap A_2$; see Figure 3.2. 
In the case that $L$ is as Figure 3.1(b), we can deform 
$L$ into the link as in Figure 3.2(a) by two $\Gamma$-moves, 
one is applied to $L\cap A_1$ and the other to $L\cap A_2$; 
see Figure 3.3. It follows from this and Figure 3.2 that 
$L$ can be deformed into a trivial link by $\Gamma$-moves, 
one is applied to $L\cap A_1$ and the others to $L\cap A_2$. 

Suppose that $p_1$ and $p_2$ are type $(i)$. 
Let $\alpha$ be an arc in the interior of $A_i$ that connects 
two singular points $p_1$ and $p_2$ of type $(i)$, and 
let $(S^3,L)=(\partial N(\alpha),
\partial(N(\alpha)\cap A_i))$. 
By the argument similar to that in the above, 
$L$ can be deformed into a trivial link by applying $\Gamma$-moves 
to $L\cap A_i$. 

Therefore, by Lemma 2.3, we can constract pass-annuli in $S^3\times[0,1]$ 
between $l$ and $l'$. Lemma 2.1 completes the proof. $\Box$

\medskip
\begin{center} 
\includegraphics[trim=0mm 0mm 0mm 0mm, width=.8\linewidth]
{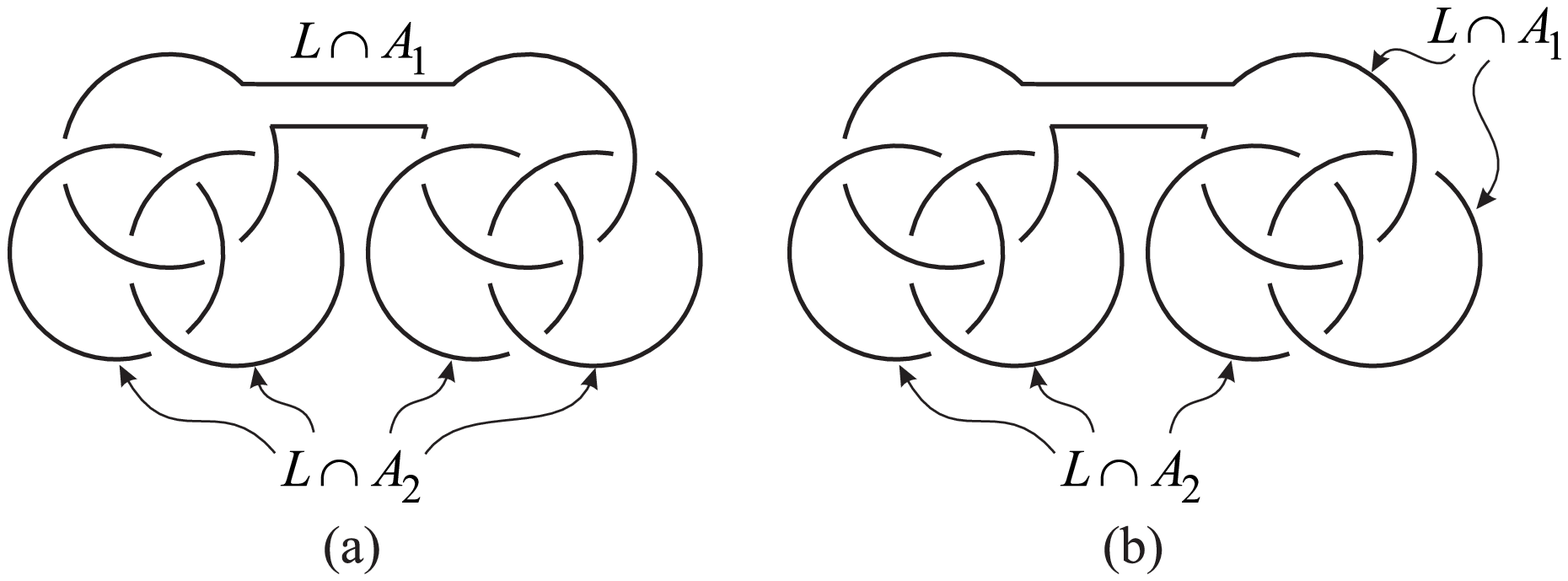}

Figure 3.1
\end{center}

\medskip
\begin{center} 
\includegraphics[trim=0mm 0mm 0mm 0mm, width=.8\linewidth]
{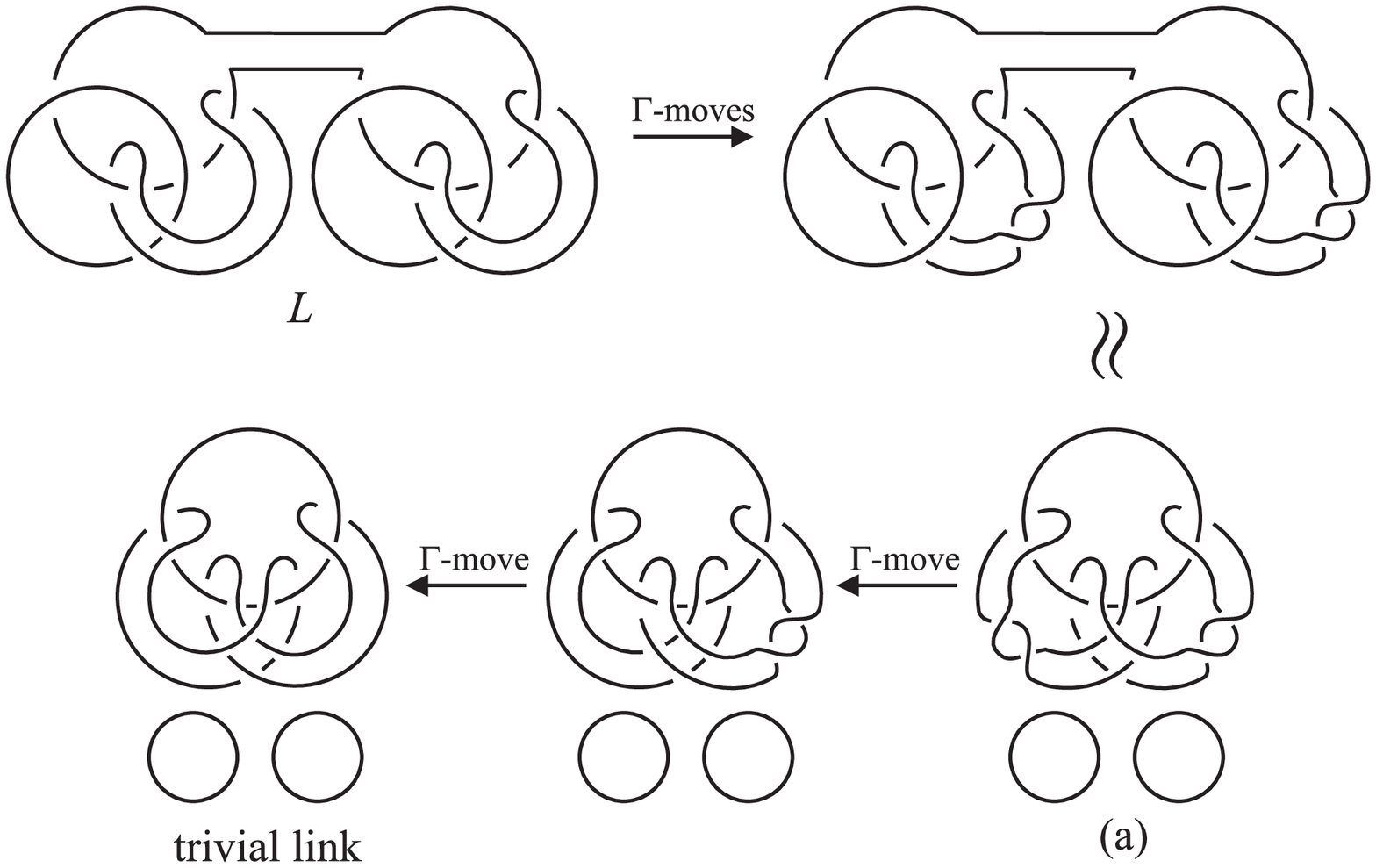}

Figure 3.2
\end{center}

\medskip
\begin{center} 
\includegraphics[trim=0mm 0mm 0mm 0mm, width=.8\linewidth]
{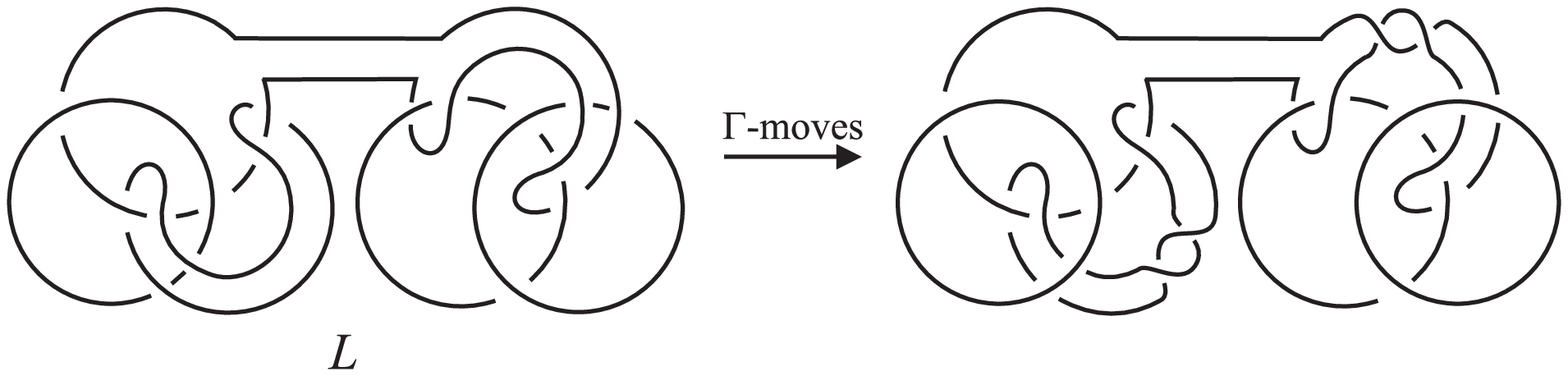}

Figure 3.3
\end{center}


\medskip
{\bf Proof of Theorem 1.1.} Since both self pass-move and 
self $\#$-move realized by link-homotopy, 
the \lq only if' parts of (i) and (ii) 
follow from \cite[Proposition]{Shi3}. 
We shall prove the \lq if' parts. 

(i) For a link $l=k_1\cup\cdots k_n$, let $G_l^o$ 
(resp. $G_l^e$) be a graph with the vertex set $\{k_1,...,k_n\}$ 
and the edge set $\{k_ik_j|\lk(k_i,k_j)\mbox{ is odd}\}$
(resp. $\{k_ik_j|\lk(k_i,k_j)\mbox{ is even}\}$). 
Note that $G_l^o\cup G_l^e$ is the complete graph with $n$ vertices. 
For a band sum $L=K_1\cup\cdots\cup K_n
(=(k_1\#_{b_1}(-\overline{k'}_1))\cup\cdots\cup
(k_n\#_{b_n}(-\overline{k'}_n)))$ of $l$ and $-\overline{l'}$, 
let $A_L$ be a graph with  the vertex set $\{K_1,...,K_n\}$ and 
the edge set $\{K_iK_j|\arf(K_i\cup K_j)=0\}$. 
(Note that $L$ is a ${\Bbb Z}_2$-algebraically split link since 
$l$ and $l'$ are link-homotopic.)

\medskip
{\bf Claim.} There is a band sum $L$ of $l$ and $-\overline{l'}$ 
such that $L$ is link-homotopic to a trivial link and $A_L$ is the 
complete graph with $n$ vertices.

\medskip
{\bf Proof.}
Let $T$ be a maximal subgraph of $G_l^o$ that does not contain 
a cycle. Since $T$ does not contain a cycle, by Lemmas 2.4 and
2.6, there is a band sum $L$ of $l$ and $l'$ such that 
$L$ is link-homotopic to a trivial link and 
$T\subset h(A_L)$, where $h:A_L\longrightarrow 
G_l^o\cup G_l^e$ the natural map defined by $h(K_i)=k_i$ and 
$h(K_iK_j)=k_ik_j$. By Lemma 2.5, we have $G_l^e\subset h(A_L)$. 
Since $h$ is injective and 
$G_l^o\cup G_l^e$ is the complete graph, it is sufficient to 
prove that $h$ is surjective. 
Let $E$ be the set of edges which are not contained in $h(A_L)$, 
and $H^o=h(A_L)\cap G_l^o$. Suppose $E\neq\emptyset$. Then 
there is an edge $e\in E$ such that there is a cycle $C$ in 
$H^o\cup e$ containning $e$ whose any \lq diagonals' are not contained 
in $G_l^o$. (In fact, for each $e_i\in E$, consider the minimum
length $l_i$ of cycles in $H^o\cup e_i$ containning $e_i$ and choose 
an edge $e$ and a cycle $C$ in $H^o\cup e$ containning $e$ 
so that length $C$ is equal to 
$\min\{l_i|e_i\in E\}$.) Without loss of generality we 
 may assume that 
$C=k_1k_2...k_ck_1$ and $e=k_1k_2$. Set $l_c=k_1\cup\cdots\cup k_c$ 
and $L_c=K_1\cup\cdots\cup K_c$. Since $C$ has no diagonals in $G_l^o$, 
all diagonals are in $G_l^e$. Thus we have 
$k_ik_j\subset H^o\cup G_l^e(=h(A_L))$ for any $i,j$ $(1\leq i<j\leq c)$ 
except for $(i,j)=(1,2)$. This implies $\arf(K_i\cup K_j)=0$ 
for any $i,j$ $(1\leq i<j\leq c,\ (i,j)\neq(1,2))$. 
The fact that $C$ has no diagonals in $G_l^o$ implies 
$l_c$ is a propre link. By the hypothesis about the Arf invariants
and Lemma 2.5, we have 
$\arf(L_c)\equiv 2\arf(l_c)\equiv 0$ (mod 2) and 
$\arf(K_i)\equiv 2\arf(k_i)\equiv 0$ (mod 2) $(i=1,...,c)$. 
Since $L_c$ is link-homotopic to a 
trivial link, by Theorem 1.2, $\arf(K_1\cup K_2)=0$. 
This contradicts $e=k_1k_2\in E$. $\Box$

\medskip
By Claim, there is a band sum $L=K_1\cup\cdots\cup K_n$ of 
$l$ and $-\overline{l'}$ such that $L$ is link-homotopic to 
a trivial link, $\arf(K_i)=0$ $(i=1,...,n)$ and $\arf(K_i\cup K_j)=0$ 
$(1\leq i<j\leq n)$. By Lemma 3.1, $L$ is self pass-equivalent to 
a trivial link. Since $L$ is a band sum of $l$ and $-\overline{l'}$, 
we can constract a pass-annuli between $l$ and $l'$. 
Lemma 2.1 completes the proof. 

(ii) Since a $\#$-move changes the value of the Arf invariant \cite{Mur}, 
by applying 
self $\#$-moves, we may assume that $\arf(k_i)=\arf(k'_i)$ for any $i$. 
Theorem 1.1(i) and Lemma 2.2 complete the proof. $\Box$

\bigskip
\footnotesize{
 }
\end{document}